\documentclass[11pt]{amsart}
\usepackage{amsmath}
\numberwithin{equation}{section}
\newtheorem{theorem}{Theorem}

\newtheorem{corollary}[theorem]{Corollary}
\newtheorem{proposition}[theorem]{Proposition}
\newtheorem{definition}{Definition}{\rm}

\def\e{\hbox{\rm e}}
\def\co{\hbox{\rm Co}}
\def\eq{\hbox{\rm Ec}}
\def\C{\mathbb{C}}
\def\J{\mathbb{J}}
\def\N{\mathbb{N}}
\def\R{\mathbb{R}}
\def\Z{\mathbb{Z}}
\def\G{\mathcal{G}}
\def\om{\mathbf{\Omega}}
\begin{document}

\title[Explicit formula for counting lattice points of polyhedra]{Simple explicit formula for counting lattice points of polyhedra}
\author{Jean B. Lasserre}
\address{LAAS-CNRS and Institute of Mathematics, 
LAAS 7 Avenue du Colonel Roche,
31077 Toulouse C\'{e}dex 4, France.}
\email{lasserre@laas.fr}

\author{Eduardo S. Zeron}
\address{Depto. Matem\'aticas, CIVESTAV-IPN,
Apdo.~Postal 14740, Mexico D.F. 07000, M\'exico.}
\email{eszeron@math.cinvestav.mx}

\date{}

\begin{abstract}
Given $z\in\C^n$ and $A\in\Z^{m{\times}n}$, we
consider the problem of evaluating the counting function
$h(y;z):=\sum\{\,z^x\,|\,x{\in}\Z^n;Ax{=}y,x{\geq}0\}$. We provide an explicit
expression for $h(y;z)$ as well as an algorithm with possibly numerous
but very simple calculations. In addition, we exhibit {\it finitely many} fixed
convex cones of $\R^n$ explicitly and exclusively defined by $A$ such
that for {\it any} $y\in\Z^m$, the sum $h(y;z)$ can be obtained by a simple
formula involving the evaluation of $\sum z^x$ over the integral points of
those cones only. At last, we also provide an alternative (and different)
formula from a decomposition of the generating function into simpler
rational fractions, easy to invert.\\
{\bf Keywords:} Computational geometry; lattice polytopes.
\end{abstract}

\maketitle

\section{introduction}
Consider the (not necessarily compact) polyhedron
\begin{equation}\label{eq1}
\om(y)\,:=\,\{x\in\R^n\;|\;Ax=y;\;x\geq0\},
\end{equation}
with $y\in\Z^m$ and $A\in\Z^{m\times{n}}$ of maximal rank for
$n\geq{m}$; besides, given $z\in\C^n$,  let $h:\Z^m\to\C$ be the
\textit{counting} function
\begin{equation}\label{eq2}
y\,\mapsto\,h(y;z)\,:=\sum_{x\in\om(y)\cap\Z^n}z^x
\end{equation}
(where $z^x$ stands for $\prod_k{z_k}^{x_k}$).
The complex vector $z\in\C^n$ may be chosen close enough to zero in
order to ensure that $h(y;z)$ is well defined even when $\om(y)$
is not compact.
If $\om(y)$ is compact, then $y\mapsto{h(y;z)}$ 
provides us with the exact number of points in the set $\om(y)\cap\Z^n$
by either evaluating $h(y,1)$, or even rounding $h(y;z)$ up to the nearest
integer when all the entries of $z$ are close enough to one.

Computation of $h$ has attracted a lot of attention in recent
years, from both theoretical and practical computation viewpoints.
Barvinok and Pommersheim \cite{barvinok2},
Brion and Vergne \cite{brion3}, 
have provided nice exact (theoretical) formulas for $h(y;z)$;
see also Szenes and Vergne \cite{szenes1}. For instance,
Barvinok considers $z\mapsto h(y;z)$ as the generating
function (evaluated at $z:=\e ^{c}\in\C^n$) of the indicator function
$x\mapsto I_{\om(y)\cap\Z^n}(x)$ of the set $\om(y)\cap\Z^n$ and
provides a decomposition into a sum of simpler generating functions
associated with supporting cones (themselves having a signed
decomposition into unimodular
cones).  We call this a {\it primal} approach because $y$ is {\it fixed},
and one works in the primal space $\R^n$ in which $\om(y)$ is defined.
Remarkably, Barvinok's counting algorithm
which is implemented in the software {\tt LattE} (see De Loera et al.
\cite{latte})
runs in time polynomial in the problem size when the dimension $n$ is
fixed. The software developed by Verdoolaege \cite{skimo} extends the
{\tt LattE} software to handle {\it parametric polytopes}. On the other hand,
Brion and Vergne \cite{brion3} consider the generating function
$H:\C ^m\to\C$ of $y\mapsto h(y;z)$, that is,
\begin{equation}
\label{gen}
w\,\mapsto\,H(w)\,:=\,\sum_{y\in\Z^m}h(y;z)w^y\,=\,\prod_{k=1}^n\frac{1}{1-z_kw^{A_k}}.
\end{equation}
They provide a generalized residue formula, and so obtain
$h(y;z)$ in closed form by {\it inversion}. We call this latter approach
{\it dual} because $z$ is fixed, and one works in the space
$\C^m$ of variables $w$ associated with the $m$ constraints $Ax=y$.

As a result of both primal and dual approaches, $h(y;z)$ is finally
expressed as a weighted sum over the vertices of $\om(y)$.
Similarly, Beck \cite{beck}, and Beck, Diaz and Robins \cite{beck2}
provided a complete analysis based on residue techniques for the
case of a tetrahedron ($m=1$).
Despite its theoretical interest, Brion and Vergne's formula is not
directly {\it tractable} because it contains many products with complex
coefficients (roots of unity) which makes the formula difficult to
evaluate numerically.  However, in some cases, this formula can be
exploited to yield an efficient algorithm as e.g. in \cite{baldoni2} for
flow polytopes, in \cite{beck3} for transportation polytopes, and more
generally when the matrix $A$ is totally unimodular as in \cite{cochet}.
Finally, in \cite{lasserre1,lasserre2}, we have provided two algorithms
based on Cauchy residue techniques to invert  $H$ in (\ref{gen}), and
an alternative algebraic technique based on partial fraction expansion
of $H$. A nice feature of the latter technique of \cite{lasserre2} is to
avoid  computing residues.
\vspace{0.2cm}

\noindent
{\bf Contribution:} Our contribution is twofold as it is concerned with
both primal and dual approaches. On the primal side, we provide an
explicit expression of $h(y;z)$ and an algorithm which involves only
elementary operations. It uses Brion's identity along with an explicit
description of the supporting cones at the vertices of $\om(y)$.
It also has a simple equivalent formulation as a (finite) {\it group
problem}. Finally, we exhibit {\it finitely many} fixed convex {\it cones}
of $\R^n$, explicitly and exclusively defined from $A$,
such that for {\it any} $y\in\Z^m$, the sum $h(y;z)$ is obtained by
a simple formula which evaluates $\sum z^x$ over the
integral points of those cones only.

On the dual side, we analyze the \textit{counting} function $h$, via
its generating function $H$ in (\ref{gen}). Inverting $H$ is difficult in
general, except if an appropriate expansion of $H$ into
simple fractions is available, as in e.g. \cite{lasserre2}. 
In their landmark paper
\cite{brion3}, Brion and Vergne provided a {\it generalized residue}
formula which yields the generic expansion
\begin{equation}\label{eq4d}
H(w)\,=\,\sum_{\sigma\in\J_A}\;\sum_{g\in{G_\sigma}}
\widehat{Q}_{g,\sigma}\prod_{k\in\sigma}
\frac{[z_kw^{A_k}]^{\delta_{k,\sigma}}}
{1-\rho_q^{g_k}[z_kw^{A_k}]^{1/q}}.
\end{equation}
Here, $\sigma\in\J_A$ whenever $A_\sigma$ is invertible,
$q$ is the smallest common multiple of all $|\det{A_\sigma}|\neq0$,
$\rho_q=\hbox{e}^{2\pi{i}/q}$ is the $q$-root of
unity, $\delta_{k,\sigma}\in\{0,1/q\}$, and $\widehat{Q}_{g,\sigma}\in\C$.
The finite group $G_\sigma$
has $q^m$ elements. The coefficients $\widehat{Q}_{g,\sigma}$ are
difficult to evaluate. Our contribution is to expand $H$ in (\ref{gen})
in the form
\begin{equation}
\label{eq4e}
H(w)\,=\,\sum_{\sigma\in\J_A}
\Bigg[\prod_{j\in\sigma}\frac{1}{1-z_jw^{A_j}}\Bigg]\times
\frac{1}{R_2(\sigma;z)}\,
\sum_{u_{\not\sigma}\in\Z_{\mu_\sigma}^{n-m}}
z^{\eta[\sigma,u_{\not\sigma}]}\,w^{A\eta[\sigma,u_{\not\sigma}]},
\end{equation}
where: $\Z_{\mu_\sigma}=\{0,1,\ldots,\mu_\sigma-1\}$,
$\mu_\sigma=|\det{}A_\sigma|$,
each $\eta[\sigma,u_{\not\sigma}]\in\Z^n$ and:
\begin{equation}\label{eqn66}
z\,\mapsto\,R_2(\sigma;z)\,:=\,\prod_{k\notin\sigma}\left[1-
\bigl(z_kz_\sigma^{-A_\sigma^{-1}\!A_k}\bigr)^{\mu_\sigma}\right].
\end{equation}
Identity (\ref{eq4e}) is a nontrivial simplification of the residue formula
(\ref{eq4d}) because the $\eta[\sigma,u_{\not\sigma}]$'s are given {\it explicitly}. And so the coefficients of the rational fraction (\ref{eq4e}) in
$w$ are very simple to evaluate with no root of unity involved (it can
also be done symbolically); however this task can be tedious as for
each $\sigma\in\J_A$ one has $|\det{}A_\sigma|^{n-m}$ terms
$\eta[\sigma,u_{\not\sigma}]$ to determine. But once determined,
(\ref{eq4e}) is easy to invert and provides $h(y;z)$ for {\it any} $y\in\Z^m$.

\section{Brion's decomposition}
\subsection{Notation and definitions} 

The notation $\C$, $\R$ and $\Z$ stand for the usual sets of complex, 
real and integer numbers, respectively. Moreover, the set of natural 
numbers $\{0,1,2,\ldots\}$ is denoted by $\N$, and for every natural 
number $\mu\in\N$, the finite set $\{0,1,\ldots,\mu-1\}$ of cardinality 
$\mu$ is denoted by $\Z_\mu$. The notation $B'$ stands for the transpose 
of a matrix (or vector) $B\in\R^{s\times{t}}$; and the $k$th column of 
the matrix $B$ is denoted by $B_k:=(B_{1,k},\ldots,B_{s,k})'$.
When $y=0$, the cone $\om(0)$ in (\ref{eq1}) is 
convex, and its \textit{dual} cone is given by,
\begin{equation}\label{eq5}
\om(0)^*\,:=\,\{b\in\R^n\,|\,b'x\geq0\;
\hbox{for every}\;x\in\om(0)\}.
\end{equation}
Notice that $\om(0)^*\equiv\R^n$ if $\om(0)=\{0\}$,
which is the case if $\om(y)$ is compact.

\begin{definition}
\label{def1}
{\rm Let $A\in\Z^{m\times{n}}$ be of maximal rank. An ordered set 
$\sigma=\{\sigma_1,\ldots,$ $\sigma_m\}$ of natural numbers is said 
to be a \textit{basis} if it has cardinality $|\sigma|=m$, the sequence of 
inequalities $1\leq\sigma_1<\sigma_2<\cdots<\sigma_m\leq{n}$ holds, and 
the square $[m\times{m}]$ submatrix~:
\begin{equation}\label{eq6} 
A_\sigma\,:=\,[A_{\sigma_1}|A_{\sigma_2}|\cdots|A_{\sigma_m}] 
\quad\hbox{is invertible}. 
\end{equation} 
We denote the set of all bases $\sigma$ by $\J_A$.}
\end{definition}

\begin{definition} 
\label{def2}
{\rm Given a maximal rank matrix $A\in\Z^{m\times{n}}$, and any basis
$\sigma\in\J_A$, the complementary matrices $A_\sigma\in\Z^{m\times{n}}$ 
and $A_{\not\sigma}\in\Z^{m\times(n-m)}$ stand for 
$[A_k]_{k\in\sigma}$ and $[A_k]_{k\notin\sigma}$, respectively. 
Similarly, given $z\in\C^n$, the complementary vectors 
$z_\sigma\in\C^m$ and $z_{\not\sigma}\in\C^{n-m}$ stand for 
$(z_k)_{k\in\sigma}$ and $(z_k)_{k\notin\sigma}$, respectively.}
\end{definition}

For each basis $\sigma\in\J_A$ with associated matrix 
$A_\sigma\in\Z^{m\times{m}}$, introduce the \textit{indicator} 
function $\delta_\sigma:\Z^m\to\N$ defined by~:
\begin{equation}\label{eq7}
y\,\mapsto\,\delta_\sigma(y)\,:=\,\left\{\begin{array}{cl}
1&\hbox{if}\;A_\sigma^{-1}y\in\Z^m,\\ 
0&\hbox{otherwise}.\end{array}\right.
\end{equation} 

Notice that $\delta_\sigma$ is a multi-periodic function with periods 
$A_\sigma$ and $\mu_\sigma:=|\det{A_\sigma}|$, meaning that 
$\delta_\sigma(y+A_{\sigma}q)=\delta_\sigma(y+\mu_{\sigma}q)=\delta_\sigma(y)$ for all $y,q\in\Z^m$. Finally, 
given a triplet $(z,x,u)\in\C^n\times\Z^n\times\R^n$, introduce the notation~:
\begin{equation}\label{eq8}
\begin{array}{rcl}
z^x&:=&z_1^{x_1}\,z_2^{x_2}\cdots{z_s}^{x_n},\\
\|z\|&:=&\max\,\{|z_1|,|z_2|,\ldots,|z_n|\},\\
\ln\langle{z}\rangle&:=&(\ln(z_1),\ln(z_2),\ldots,\ln(z_n)).
\end{array}
\end{equation}

Notice that $z^x=z_{\sigma}^{x_{\sigma}}z_{\not\sigma}^{x_{\not\sigma}}$, for all bases 
$\sigma\in\J_A$ and all $z\in\C^n,\,x\in\Z^n$.

\subsection{Brion's decomposition} 

Let $\om(y)$ be the convex polyhedron in (\ref{eq1}) with
$y\in\Z^m,A\in\Z^{m\times{n}}$ being of maximal rank, and let  
$h:\Z^m\to\C$ be the counting function in (\ref{eq2}), 
with $\|z\|<1$.

Obviously $h(y;z)=0$ whenever the equation $Ax=y$ has no 
solution $x\in\N^n$. The main idea is to decompose the function 
$h$ following Brion's ideas. Given any convex 
rational polyhedron $P\subset\R^n$, let $[P]:\R^n\to\{0,1\}$ be 
its characteristic function, and $f[P]:\C\to\C$ its associated 
rational function, such that
\begin{equation}\label{eq9}
z\,\mapsto\,f[P,z]\,:=\sum_{x\in{P}\cap\Z^n}z^x,
\end{equation}
holds whenever the sum converges absolutely. For every vertex $V$ of 
$P$, define $\co(P,V)\subset\R^n$ to be the supporting cone of $P$ at 
$V$. Then, Brion's formula yields the decomposition~:
\begin{equation}\label{eq10-epsilon}
[P]\,=\,\sum_{\hbox{\footnotesize vertices}\:V}[\co(P,V)],
\end{equation}
modulo the group generated by the characteristic functions 
of convex polyhedra which contain affine lines. And so,
\begin{equation}\label{eq10}
f[P,z]\,=\,\sum_{\hbox{\footnotesize vertices}\:V}f[\co(P,V),z].
\end{equation}

The above summation is \textit{formal} because
in general there is no $z\in\C^n$ for which the series 
$$\sum\{z^x\,|\,x\in{P}\cap\Z^n\}\quad\hbox{and}
\quad\sum\{z^x\,|\,x\in\co(P,V)\cap\Z^n\}$$
converge absolutely for all vertices $V$.
The notation $\sum{E}$ stands for the sum of all elements of a 
countable set $E\subset\C$. It is a complex number 
whenever the resulting series converges absolutely; otherwise it 
stands for a formal series. 

\vspace{9pt}\noindent
\textbf{Example:} Let $P:=[0,1]\subset\R$ so that
$\co(P,\{0\})=[0,+\infty)$ and $\co(P,\{1\})=(-\infty,1]$.
Simple enumeration yields $f[P,z]=z^0+z=1+z$, but one also has:
\[f[P,z]\,=\,f[\co(P,\{0\}),z]+f[\rm(P,\{1\}),z]\,=\,1/(1-z)+z^2/(z-1)\,=\,1+z.\]

\section{Computing $h(y;z)$: A primal approach}
\label{section-primal}

Let $C(\J_A):=\{Ax\,|\,x\in\N^n\}\subset\R^m$ be the cone generated by the 
columns of $A$, and for any basis $\sigma\in\J_A$, let $C(\sigma)\subset\R^m$ 
be the cone generated by the columns $A_k$ with $k\in\sigma$. As $A$ 
has maximal rank, $C(\J_A)$ is the union of all 
$C(\sigma)$, $\sigma\in\J_A$. With any $y\in{}C(\J_A)$ associate 
the intersection of all cones $C(\sigma)$ that contain $y$. This defines a subdivision of $C(\J_A)$ into polyhedral cones. 
The interiors of the maximal subdivisions are called {\it chambers}. 
In each chamber $\gamma$, the polyhedron $\om(y)$ is {\it simple}, i.e.  
$A_\sigma^{-1}y>0$ for all $\sigma\in\J_A$ such  that $A_\sigma^{-1}y\geq0$.

For any chamber $\gamma$, define,
\begin{equation}\label{supp-base}
\mathcal{B}(\J_A,\gamma)\,:=\,\{\sigma\in\J_A\:|\:\gamma\subset{}C(\sigma)\}.
\end{equation}
The intersection of all $C(\sigma)$ with $\sigma\in\mathcal{B}(\J_A,\gamma)$
is the closure $\overline{\gamma}$ of $\gamma$.

Back to our original problem, and setting $P:=\om(y)$, the rational 
function $f[P,z]$ is equal to $h(y;z)$ in (\ref{eq2}) whenever $\|z\|<1$. We 
next provide an explicit description of the rational function $f[\co(P,V),z]$ 
for every vertex $V$ of $P$.

Let $\delta_\sigma$ be the function defined in (\ref{eq7}), 
and let $\Z_{\mu_\sigma}:=\{0,1,\ldots,\mu_\sigma-1\}$ with 
$\mu_\sigma:=|\det{A_\sigma}|$.  A vector $V\in\R^n$ is a vertex of 
$P=\om(y)$ if and only if there exists a basis $\sigma\in\J_A$ such that~:
\begin{equation}\label{supp-vertex}
V_\sigma\,=\,A_\sigma^{-1}y\,\geq\,0\quad\hbox{and}\quad 
V_{\not\sigma}=0,
\end{equation}
where $V_\sigma$ and $V_{\not\sigma}$ are given in Definition~\ref{def2}.
Moreover, the supporting cone of $P$ at the vertex $V$ is described by~:
\begin{equation}\label{supp-cone}
\co(\om(y),V)\,:=\,\left\{x\in\R^n\:|\:
Ax=y;\;x_k\geq0\;\mbox{if}\;V_k=0\right\}.
\end{equation}
Let us now define the larger set
\begin{equation}\label{eq14} 
C(\om(y),\sigma)\,:=\,\{x\in\R^n\,|\,A_{\sigma}x_{\sigma}+
A_{\not\sigma}x_{\not\sigma}=y;\;x_{\not\sigma}\geq0\},
\end{equation}
so that $\co(\om(y),V)$ is a subcone of $C(\om(y),\sigma)$ for 
all bases $\sigma\in\J_A$ and vertex $V$ of $\om(y)$ which satisfy 
$V_{\not\sigma}=0$ (recall (\ref{supp-vertex})). 
Besides, when $V_{\not\sigma}=0$ and $y\in\gamma$ for some chamber 
$\gamma$, then $C(\om(y),\sigma)$ and $\co(\om(y),V)$ 
are identical because $\om(y)$ is a simple polytope, and so 
$A_\sigma^{-1}y>0$ for all $\sigma\in \J_A$.

Recall that $A_\sigma\in\Z^{m\times{n}}$ and 
$A_{\not\sigma}\in\Z^{m\times(n-m)}$ stand for 
$[A_k]_{k\in\sigma}$ and $[A_k]_{k\notin\sigma}$,  respectively. 
Similarly, given a vector $x\in\Z^n$, the vectors 
$x_\sigma$ and $x_{\not\sigma}$ 
stand for $(x_k)_{k\in\sigma}$ and 
$(x_k)_{k\notin\sigma}$ respectively.
The following result is from \cite[p.~818]{brion3}.

\begin{proposition}
\label{prop-cone}
Let $y\in\R^m$ and let $\om(y)$ be as in (\ref{eq1}), and 
let $y\in\overline{\gamma}$ for some chamber $\gamma$. Then,
\begin{equation}\label{prop-cone-1}
[\om(y)]\,=\,\sum_{\sigma\in\mathcal{B}(\J_A,\gamma)}
[C(\om(y),\sigma)],
\end{equation}
modulo the group generated by the characteristic 
functions of convex polyhedra which contain affine lines. 
\end{proposition}

\begin{proof}
Using notation of \cite[p.~817]{brion3}, define the linear 
mapping $p:\R^n\to\R^m$ with $p(x)=Ax$, so that the polyhedra 
$P_{\Delta}(y)$ and $\om(y)$ are identical. Moreover, for every 
basis $\sigma\in\mathcal{B}(\J_A,\gamma)$, 
$v_\sigma:\R^m\to\R^n$ is the linear 
mapping:
$$y\mapsto\quad [v_\sigma(y)]_\sigma\,=\,A_\sigma^{-1}y\quad
\hbox{and}\quad[v_\sigma(y)]_{\not\sigma}\,=\,0,\qquad y\in\R^m.$$
  
Finally, for every $\hat{x}\in\R^n$ with $\hat{x}\geq0$, 
$\rho_{\sigma}(\hat{x}):=\hat{x}-v_{\sigma}(A\hat{x})$
satisfies, 
$$[\rho_\sigma(\hat{x})]_\sigma=-A_\sigma^{-1}
A_{\not\sigma}\hat{x}_{\not\sigma}\quad\hbox{and}\quad
[\rho_\sigma(\hat{x})]_{\not\sigma}\,=\,\hat{x}_{\not\sigma}.$$

Therefore, the cone $[v_\sigma(y)+\rho_\sigma(C)]$ in \cite{brion3} is 
the set of points $x\in\R^m$ such that $x_{\not\sigma}\geq0$ and 
$x_{\sigma}=A_\sigma^{-1}(y-A_{\not\sigma}x_{\not\sigma})$; and 
so this cone is just $[C(\om(y),\sigma)]$ in~(\ref{eq14}).
Therefore a direct application of $(3.2.1)$ in \cite[p.~818]{brion3} 
yields (\ref{prop-cone-1}).
\end{proof}
\begin{theorem}\label{Brion}
Let $y\in\Z^m,\,z\in\C^n$ with $\|z\|<1$, and let 
$y\in\overline{\gamma}$ for some chamber $\gamma$. Recall the set 
of bases $\mathcal{B}(\J_A,\gamma)$ defined in (\ref{supp-base}). 
With $P:=\om(y)$, the rational function $h$ defined in (\ref{eq2}) 
can be written:
\begin{equation}\label{eq15} 
h(y;z)\,=\sum_{\sigma\in\mathcal{B}(\J_A,\gamma)}
f[C(\om(y),\sigma),z]\,=\sum_{\sigma\in\mathcal{B}
(\J_A,\gamma)}\;\frac{R_1(y,\sigma;z)}{R_2(\sigma;z)},
\end{equation}
\begin{eqnarray}\label{eq11}
with\qquad z\,\mapsto\,R_1(y,\sigma;z)&:=&z_\sigma^{A_\sigma^{-1}y}
\sum_{u\in\Z_{\mu_\sigma}^{n-m}}
\frac{\delta_\sigma(y-A_{\not\sigma}u)\,z_{\not\sigma}^u}
{z_\sigma^{A_\sigma^{-1}\!A_{\not\sigma}u}},\\
\label{eq12}
and\qquad z\,\mapsto\,R_2(\sigma;z)&:=&\prod_{k\notin\sigma}\left[1-\bigl(
z_kz_\sigma^{-A_\sigma^{-1}\!A_k}\bigr)^{\mu_\sigma}\right].
\end{eqnarray}
\end{theorem}
The pair $\{R_1,R_2\}$ is well defined whenever $z\in\C^n$ 
satisfies $z_k\neq0$ and $z_k\neq{z_\sigma}^{A_\sigma^{-1}\!A_k}$ for 
every basis $\sigma\in\J_A$ which does not contain the index
$k\not\in\sigma$.

\begin{proof} By a direct application of Brion's theorem to the sum 
(\ref{prop-cone-1}), the associated rational functions  
$f[\om(y),z]$ and $f[C(\om(y),\sigma),z]$ satisfy:
\begin{equation}\label{eq16-epsilon}
h(y,z)\,=\,f[\Omega(y),z]\,=\sum_{\sigma\in
\mathcal{B}(\J_A,\gamma)}f[C(\om(y),\sigma),z].
\end{equation}

Therefore, in order to show (\ref{eq15}), one only needs to prove 
that the rational function $\frac{R_1(y,\sigma;z)}{R_2(\sigma;z)}$ 
is equal to $f[C(\om(y),\sigma),z]$, i.e., 
\begin{equation}\label{eq16}
\frac{R_1(y,\sigma;z)}{R_2(\sigma;z)}\,=\,
\sum\,\{z^x\,|\,x\in C(\Omega(y),\sigma)\cap\Z^n\},
\end{equation}
on the domain 
$D_\sigma\,=\,\{z\in\C^n\,|\,1>|z_kz_\sigma^{-A_\sigma^{-1}\!A_k}|\;
\hbox{for each}\;k\not\in\sigma\}$.
Notice that
\begin{eqnarray*}
&&\frac{1}{R_2(\sigma;z)}\;=\;\prod_{k\notin\sigma}\;\frac{1}
{1-\bigl(z_kz_\sigma^{-A_\sigma^{-1}\!A_k}\bigr)^{\mu_\sigma}}\,=\\
&&=\;\prod_{k\notin\sigma}\;\sum_{v_k\in\N}\biggl[\frac{z_k}
{z_\sigma^{A_\sigma^{-1}\!A_k}}\biggr]^{\mu_{\sigma}v_k}\;=\;
\sum_{v\in\N^{n-m}}\;\frac{z_{\not\sigma}^{\mu_{\sigma}v}}
{z_\sigma^{\mu_{\sigma}A_\sigma^{-1}\!A_{\not\sigma}v}},
\end{eqnarray*}
on $D_\sigma$. On the other hand,
according to~(\ref{eq14}), the integer vector $x\in\Z^n$ lies 
inside the cone $C(P,V_\sigma)$ if and only if~:
\begin{eqnarray*}
x_\sigma&=&A_\sigma^{-1}(y-A_{\not\sigma}x_{\not\sigma}),\quad
\delta_\sigma(y-A_{\not\sigma}x_{\not\sigma})\;=\;1\quad\hbox{and}\\
\nonumber x_{\not\sigma}&=&u+\mu_{\sigma}v,\quad\hbox{with}\quad
u\in\Z_{\mu_\sigma}^{n-m}\quad\hbox{and}\quad v\in\N^{n-m}.
\end{eqnarray*}
But from the definition~(\ref{eq11}) of $R_1(y,\sigma;z)$ and 
$z^x=z_{\not\sigma}^{x_{\not\sigma}}z_{\sigma}^{x_{\sigma}}=z_{\sigma}^{A_\sigma^{-1}(y-A_{\not\sigma}x_{\not\sigma})}$,
\begin{eqnarray}\label{R1/R2}
&&\frac{R_1(y,\sigma;z)}{R_2(\sigma;z)}\;=\;z_\sigma^{A_\sigma^{-1}y}
\sum_{u\in\Z_{\mu_\sigma}^{n-m}}\;\sum_{v\in\N^{n-m}}\frac{
\delta_\sigma(y-A_{\not\sigma}u)\,z_{\not\sigma}^{x_{\not\sigma}}}
{z_\sigma^{A_\sigma^{-1}\!A_{\not\sigma}x_{\not\sigma}}},\\
\nonumber&&\,=\,\sum\,\{z^x\,|\,x\in C(\om(y),\sigma)\cap\Z^n\}
\;=\;f[C(\om(y)\sigma),z],
\end{eqnarray}
which is exactly~(\ref{eq16}). Notice that $x_{\not\sigma}=u+\mu_{\sigma}v$, and so $\delta_\sigma(y-A_{\not\sigma}u)=\delta_\sigma(y-A_{\not\sigma}x_{\not\sigma})$ because of the 
definition~(\ref{eq7}) of $\delta_\sigma$. Finally, 
using (\ref{R1/R2}) in (\ref{eq16-epsilon}) yields 
that (\ref{eq15}) holds whenever $\|z\|<1$ and
$R_1(y,\sigma;z)$ and $R_2(\sigma;z)$ are all well defined.
\end{proof}

Notice that $R_2$ is constant with respect to $y$, 
and from the definition (\ref{eq7}) of $\delta_\sigma$, 
$R_1$ is \textit{quasiperiodic} with periods $A_\sigma$ 
and $\mu_\sigma$, meaning that 
\begin{equation}\label{eq13}
\begin{array}{rcl}
R_1(y+A_{\sigma}q,\sigma;z)&=&R_1(y,\sigma;z)\,
z_\sigma^q\quad\hbox{and}\\ 
R_1(y+\mu_{\sigma}q,\sigma;z)&=&R_1(y,\sigma;z)\,
\bigl(z_\sigma^{A_\sigma^{-1}q}\bigr)^{\mu_\sigma}
\end{array}\end{equation}
hold for all $y,q\in\Z^m$.
Obviously, the more expensive part in calculating 
$R_2(\cdot)$ in (\ref{eq12})  is to compute the determinant 
$\mu_\sigma=|\det{A}_\sigma|$. On the other hand, computing 
$R_1(\cdot)$ in (\ref{eq11}) may become quite expensive when
$\mu_\sigma$ is large, as one must evaluate $\mu_\sigma^{n-m}$ 
terms, the cardinality of $\Z_{\mu_\sigma}^{n-m}$. However, as 
detailed below, a more careful analysis of (\ref{eq11}) yields 
some simplifications.
\subsection{Simplifications via group theory}

From the proof of Theorem~\ref{Brion}, 
the closed forms (\ref{eq11})--(\ref{eq12}) for $R_1(\cdot)$ 
and $R_2(\cdot)$ are deduced from (\ref{R1/R2}), i.e.,
$$\frac{R_1(y,\sigma;z)}{R_2(\sigma;z)}\,=\,
z_\sigma^{A_\sigma^{-1}y}\sum_{x_{\not\sigma}\in\Z^{n-m}}
\frac{\delta_\sigma(y-A_{\not\sigma}x_{\not\sigma})\,
z_{\not\sigma}^{x_{\not\sigma}}}
{z_\sigma^{A_\sigma^{-1}\!A_{\not\sigma}x_{\not\sigma}}},$$
after setting $x_{\not\sigma}=u+\mu_{\sigma}v$ and recalling 
that $\delta_\sigma(y)$  is a periodic function, i.e.,
$\delta_\sigma(y+\mu_{\sigma}q)=\delta_\sigma(y)$ for all $y,q\in\Z^m$. However, we have not 
used yet that $\delta_\sigma(y+A_{\sigma}q)=\delta_\sigma(y)$ as well. 
For every $\sigma\in\J_A$, consider the lattice~: 
\begin{equation}
\label{lattice}
\Lambda_\sigma\,:=\,\bigoplus_{j\in\sigma}A_j\Z\,\subset\,\Z^m,
\end{equation}
generated by the columns $A_j$, $j\in\sigma$. 
The following quotient group 
\begin{eqnarray}\label{group}
\G_\sigma&:=&\Z^m/\Lambda_{\sigma}\;=\;
\Z^m\Big/\bigoplus_{j\in\sigma}A_j\Z\\
\nonumber&=&\{\eq[0,\sigma],\eq[2,\sigma],
\ldots,\eq[\mu_\sigma-1,\sigma]\}
\end{eqnarray}
is commutative, with $\mu_\sigma=|\det{A_\sigma}|$ elements 
(or, equivalence classes) $\eq[j,\sigma]$, and so, $\G_\sigma$ 
is isomorphic to a finite Cartesian product of cyclic groups 
$\Z_{\eta_k}$, i.e.,
$$\G_\sigma\,\cong\,\Z_{\eta_1}\times
\Z_{\eta_2}\times\cdots\times\Z_{\eta_s}.$$

Obviously, $\mu_\sigma=\eta_1\eta_2\cdots\eta_s$, and so, $\G_\sigma$ is isomorphic 
to the cyclic group $\Z_{\mu_\sigma}$ whenever $\mu_\sigma$ is 
a prime number. Actually, $\G_\sigma=\{0\}$ whenever $\mu_\sigma=1$. Notice that the Cartesian
product $\Z_{\eta_1}\times\cdots\times\Z_{\eta_s}$ can be seen 
as the integer space $\Z^s$ modulo the vector
$\eta\,:=\,(\eta_1,\eta_2,\cdots,\eta_s)'\,\in\N^s$.

Hence, for every finite commutative group $\G_\sigma$, 
there exist a positive integer $s_\sigma\geq1$, a vector 
$\eta_\sigma\in\N^{s_\sigma}$ with positive entries, 
and a group isomorphism,
\begin{equation}\label{eq33}
g_\sigma\,:\,\G_\sigma\,\to\,\Z^{s_\sigma}\bmod\eta_\sigma,
\end{equation}
where $g_\sigma(\xi)\bmod\eta_\sigma$ means evaluating
$[g_\sigma(\xi)]_k\bmod[\eta_\sigma]_k$, for all indices
$1\leq{k}\leq{s_\sigma}$. For every $y\in\Z^m$, 
there exists a unique equivalence class $\eq[j_y,\sigma]$ which 
contains $y$, and so we can define the following group epimorphism,
\begin{eqnarray}\label{eq34}
&&\hat{h}_\sigma\,:\,\Z^m\,\to\,\Z^{s_\sigma}\mod\eta_\sigma,\\
\nonumber&&y\,\mapsto\,\hat{h}_\sigma(y)\,:=\,g_\sigma(\eq[j_y,\sigma]).
\end{eqnarray}

On the other hand, the unit element of $\G_\sigma$
is the equivalence class $\eq[0,\sigma]$ which contains the 
origin, that is, $\eq[0,\sigma]\,=\,\{A_{\sigma}q\;|\;q\in\Z^m\}$.

Hence, $\hat{h}_\sigma(y)=0$ if and only if there
exists $q\in\Z^m$ such that $y=A_{\sigma}q$. We can then redefine 
the function $\delta_\sigma$ as follows,
\begin{equation}\label{eq35}
y\,\mapsto\,\delta_\sigma(y)\,:=\,\left\{
\begin{array}{cl}1&\hbox{if}\;\hat{h}_\sigma(y)=0,\\ 
0&\hbox{otherwise},\end{array}\right.
\end{equation}
One also needs the following 
additional notation; given any matrix $B\in\Z^{m\times{t}}$, 
\begin{equation}\label{eq36}
\hat{h}_\sigma(B)\,:=\,[\hat{h}_\sigma(B_1)|\hat{h}_\sigma(B_2)|
\cdots|\hat{h}_\sigma(B_t)]\,\in\,\Z^{s_\sigma\times{t}}.
\end{equation}

And so, from (\ref{eq11}),
$\hat{h}_\sigma(y-A_{\not\sigma}u)\,\equiv\,\hat{h}_\sigma(y)
-\hat{h}_\sigma(A_{\not\sigma})u\bmod\eta_\sigma$.
Finally, using (\ref{eq35}) in (\ref{eq11}), one obtains 
a simplified version of $R_1(\cdot)$ in the form:
\begin{equation}\label{eq37}
R_1(y,\sigma;z)\,=\,\sum\Biggl\{
\frac{z_\sigma^{A_\sigma^{-1}y}z_{\not\sigma}^u}
{z_\sigma^{A_\sigma^{-1}A_{\not\sigma}u}}\;\biggm|\;
\begin{array}{ll} u\in\Z^{n-m}_{\mu_\sigma};\\
\hat{h}_\sigma(y)\equiv\hat{h}_\sigma(A_{\not\sigma})u
\bmod\eta_\sigma\end{array}\Biggr\}.
\end{equation}

Next, with $q\in\Z^m$ fixed,
$\nu_qA_\sigma^{-1}q\in\Z^m$ for some 
integer $\nu_q$, if and only if $\nu_q\hat{h}_\sigma(q)=0\bmod
\eta_\sigma$. If we set $\nu_q=\mu_\sigma$, 
then $\mu_{\sigma}A_\sigma^{-1}q\in\Z^m$, and $\mu_\sigma\hat{h}_\sigma(q)=0\bmod\eta_\sigma$, because $\G_\sigma$ has $\mu_\sigma=|\det{A}_\sigma|$ 
elements. Nevertheless, $\mu_\sigma$ may not be the 
smallest positive integer with that property. So, given 
$\sigma\in\J_A$ and $k\notin\sigma$, define $\nu_{k,\sigma}\geq1$ 
to be \textit{order} of $\hat{h}_\sigma(A_k)$. That is, 
$\nu_{k,\sigma}$ is the smallest positive integer such that 
$\nu_{k,\sigma}\hat{h}_\sigma(A_k)=0\bmod \eta_\sigma$, or equivalently~:
\begin{equation}\label{eq40}
\nu_{k,\sigma}A_\sigma^{-1}A_k\,\in\,\Z^m.
\end{equation}
Obviously $\nu_{k,\sigma}\leq\mu_\sigma$. Moreover, $\mu_\sigma$ is a multiple 
of $\nu_{k,\sigma}$ for it is the order of an element in 
$\G_\sigma$. For example, the group $\Z^2$ modulo 
$\eta=\binom{2}{7}$ has 14 elements; and the elements 
$b_1=\binom{1}{0}$, $b_2=\binom{0}{1}$ and $b_3=\binom{1}{1}$ 
have respective orders~: 2, 7 and 14. Notice that,
$2b_1\,\equiv\,7b_2\,\equiv\,14b_3\,\equiv\,0\bmod\eta$.
But, $2b_3\,\equiv\,2b_2\,\not\equiv\,0$ and
$7b_3\,\equiv\,b_1\not\equiv\,0\bmod\eta$.

The important observation is that $\delta_\sigma(y-\nu_{k,\sigma}A_kq)=\delta_\sigma(y)$ for all $q\in\Z^m$ and $k\notin\sigma$,
which follows from (\ref{eq40}) and (\ref{eq7}). Thus, 
following step by step the proof of Theorem \ref{Brion}, we obtain:

\begin{corollary}\label{simplify}
Let $y\in\Z^m,\,z\in\C^n$ with $\|z\|<1$, and let 
$y\in\overline{\gamma}$ for some chamber $\gamma$. Recall the set 
of bases $\mathcal{B}(\J_A,\gamma)$ defined in (\ref{supp-base}). 
With $\sigma\in\mathcal{B}(\J_A,\gamma)$, let $R_1$ and $R_2$ be 
as in Theorem \ref{Brion}. Then 
\begin{equation}\label{eq41}
\frac{R_1(y,\sigma;z)}{R_2(y;z)}\,=\,
\frac{R^*_1(y,\sigma;z)}{R^*_2(y;z)},
\end{equation}
\begin{equation}\label{eq43}
where:\qquad R^*_2(\sigma;z)\,:=\,\prod_{k\notin\sigma}\left[1-\bigl(
z_kz_\sigma^{-A_\sigma^{-1}\!A_k}\bigr)^{\nu_{k,\sigma}}\right],
\end{equation}
\begin{eqnarray}\label{eq42}
&&R^*_1(y,\sigma;z)\;:=\;z_\sigma^{A_\sigma^{-1}y}
\sum_{u_{\not\sigma}\in{U}_{\not\sigma}}\frac{\delta_\sigma
(y-A_{\not\sigma}u_{\not\sigma})\,z_{\not\sigma}^{u_{\not\sigma}}}
{z_\sigma^{A_\sigma^{-1}\!A_{\not\sigma}u_{\not\sigma}}}\;=\\
\nonumber&&=\;\sum\Biggl\{
\frac{z_\sigma^{A_\sigma^{-1}y}z_{\not\sigma}^{u_{\not\sigma}}}
{z_\sigma^{A_\sigma^{-1}A_{\not\sigma}u_{\not\sigma}}}\;\biggm|\;
\begin{array}{ll} u_{\not\sigma}\in{U}_{\not\sigma};\\
\hat{h}_\sigma(y)\equiv\hat{h}_\sigma(A_{\not\sigma})u_{\not\sigma}
\mod\eta_\sigma\end{array}\Biggr\},
\end{eqnarray}
with $U_{\not\sigma}\,:=\,\{u_{\not\sigma}\in\N^{n-m}\;|\;
0\leq{u_k}\leq\nu_{k,\sigma}-1;\;k\notin\sigma\}$.
\end{corollary}

One can also obtain (\ref{eq41}) by noticing that:
$$\frac{R_1(y,\sigma;z)}{R^*_1(y,\sigma;z)}
\,=\,\frac{R_2(\sigma;z)}{R^*_2(\sigma;z)}\,=\,
\prod_{k\notin\sigma}\left(1+\beta^{\nu_{k,\sigma}}
+\cdots+\beta^{\mu_\sigma-\nu_{k,\sigma}}\right),$$
where $\beta_{k,\sigma}=z_kz_\sigma^{-A_\sigma^{-1}A_k}$,
and $\mu_\sigma$ is a multiple of $\nu_{k,\sigma}$. 

\subsection{Simplifications via finite number of generators}

Decompose $\Z^m$ into $\mu_\sigma:=|\det{A_\sigma}|$ disjoint 
equivalent classes, where $y,\xi\in\Z^m$ 
are equivalent if and only if $\delta_\sigma(y-\xi)=1$. 
For every basis $\sigma\in\J_A$, let $\G_\sigma$ be the quotient 
group defined in~(\ref{group}), that is,
$$\G_\sigma\;=\;\Z^m\Big/\bigoplus_{j\in\sigma}A_j\Z\;=\;
\{\eq[0,\sigma],\ldots,\eq[\mu_\sigma-1,\sigma]\}.$$
Notice that $y,\xi\in\Z^n$ belong to $\eq[j,\sigma]$ if and only if 
$A_\sigma^{-1}(y-\xi)\in\Z^n$, and that $\Z^m$ is equal 
to the disjoint union of all classes $\eq[j,\sigma]$. 

Next, pick up a {\it minimal} representative 
element of every class, i.e., fix 
\begin{equation}\label{eq23}
\xi[j,\sigma]\in\eq[j,\sigma]\quad\hbox{such that}\quad
A_\sigma^{-1}y\geq{A}_\sigma^{-1}\xi[j,\sigma]\geq0,
\end{equation}
for every $y\in\eq[j,\sigma]$ with $A_\sigma^{-1}y\geq0$. The 
minimal representative elements $\xi[j,\sigma]$ in~(\ref{eq23}) 
can be computed as follows: Let $d\in\eq[j,\sigma]$, arbitrary,
and let $d^*\in\Z^m$ be such that his $k$-entry 
$d^*_k$ is the smallest integer greater than or equal to the 
$k$-entry of $-A_\sigma^{-1}d$. The vector $\xi[j,\sigma]$
defined by $d+A_{\sigma}d^*$ satisfies~(\ref{eq23}).

Notice that $d^*+A_\sigma^{-1}d\geq0$. Besides, let 
$d,y\in\eq[j,\sigma]$ with $A_\sigma^{-1}y\geq0$. There exists 
$q\in\Z^m$ such that $y=d+A_{\sigma}q$. Hence $q\geq-A_\sigma^{-1}d$; 
in addition, $q\geq{d^*}$ follows from the above definition of $d^*$, 
and so $A_\sigma^{-1}y\,\geq\,d^*+A_\sigma^{-1}d\,\geq\,0$.

Therefore, the vector $\xi[j,\sigma]:=d+A_{\sigma}d^*$ 
satisfies~(\ref{eq23}). In particular, if $\eq[0,\sigma]$ is the 
class which contains the origin of $\Z^m$, then $\xi[0,\sigma]=0$. 
Notice that for every integer vector $y\in\Z^m$, there exists 
a unique $\xi[j,\sigma]$ such that~: 
$$y\,=\,\xi[j,\sigma]\,+\,A_{\sigma}\,q,
\quad\hbox{for}\quad{q}\in\Z^m.$$ 
Moreover, the extra condition $A_\sigma^{-1}y\geq0$ 
holds if and only if:
\begin{equation}\label{eq24}
y\,=\,\xi[j,\sigma]\,+\,A_{\sigma}\,q
\quad\hbox{with}\quad{q}\in\N^m.
\end{equation}

We obtain a compact form of $h(y;z)$ when 
$y\in\Z^m\cap\overline{\gamma}$, for some chamber $\gamma$.

\begin{theorem}\label{inversion}
Let $h$ and $\xi[j,\sigma]$ be as in~(\ref{eq2}) and~(\ref{eq23}), 
respectively. Let $y\in\Z^m\cap\overline{\gamma}$, for some chamber 
$\gamma$. Recall the set of bases $\mathcal{B}(\J_A,\gamma)$ defined in 
(\ref{supp-base}). For every basis $\sigma\in\mathcal{B}(\Delta,\gamma)$ 
there is a unique index  $0\leq\jmath[y,\sigma]<\mu_\sigma$ such that $y$ 
is contained in the equivalence class $\eq[\jmath[y,\sigma],\sigma]$ 
defined in~(\ref{group}), and so:
\begin{equation}\label{eq28}
h(y;z)\,=\sum_{\sigma\in\mathcal{B}(\Delta,\gamma)}
\frac{R_1(\xi[\jmath[y,\sigma],\sigma],\sigma;z)}
{R_2(\sigma;z)}\;z_\sigma^{\lfloor{}A_\sigma^{-1}y\rfloor},
\end{equation}
where $\lfloor{A_\sigma}^{-1}y\rfloor\in\Z^m$ is
such that his $k$-entry is the largest integer 
less than or equal to the $k$-entry of $A_\sigma^{-1}y$.
\end{theorem}

\begin{proof}
Recall that if $y\in\Z^m\cap\overline{\gamma}$
$$h(y;z)\,=\sum_{\sigma\in\mathcal{B}(\Delta,\gamma)}
\frac{R_1(y,\sigma;z)}{R_2(y;z)}$$
Next, recalling the definition (\ref{supp-base}) of 
$\mathcal{B}(\J_A,\gamma)$, $A_\sigma^{-1}y\geq0$ 
for every basis $\sigma\in\mathcal{B}(\J_A,\gamma)$ with 
$y\in\overline{\gamma}$. Recall that there is a unique index 
$\jmath[y,\sigma]<\mu_\sigma$ such that 
$y=\xi[\jmath[y,\sigma],\sigma]+A_\sigma{q}$ with $q\in\N^m$; 
see (\ref{eq24}) and the comment just before. 

To obtain the vector $q\in\N^m$, recall that
the minimal representative element $\xi[\jmath[y,\sigma],\sigma]$ 
in~(\ref{eq23}) is the sum $y+A_{\sigma}y^*$ 
where $y^*\in\Z^m$ is such that his $k$-entry $y^*_k$ is the smallest 
integer greater than or equal to $-A_\sigma^{-1}y$, for we only need 
to fix $d=y$ in the paragraph that follows (\ref{eq23}). In 
particular, $\lfloor{A_\sigma}^{-1}y\rfloor=-y^*$, and
$\xi[\jmath[y,\sigma],\sigma]=y-A_{\sigma}\lfloor{A_\sigma}^{-1}y\rfloor$, which when used in (\ref{eq11}) and (\ref{eq13}), yields,
$$R_1(y,\sigma;z)\,=\,R_1(\xi[\jmath[y,\sigma],\sigma],\sigma;z)
\,z_\sigma^{\lfloor{}A_\sigma^{-1}y\rfloor}.$$
And so (\ref{eq15}) implies (\ref{eq28}).
\end{proof}

Theorem \ref{inversion} explicitly shows that it suffices to compute 
$R_1(v,\sigma;z)$ for finitely many values $v=\xi[j,\sigma]$, 
with $\sigma\in\mathcal{B}(\Delta,\gamma)$ and $0\leq{j}<\mu_\sigma$, 
in order to calculate $h(y;z)$ for arbitrary values 
$y\in\Z^m\cap\overline{\gamma}$, via (\ref{eq28}). 

In other words, in the closure $\overline{\gamma}$ of a chamber 
$\gamma$, one only needs to consider {\it finitely many} fixed 
convex cones $C(\om(\xi[j,\sigma]),\sigma)\subset\R^n$, where 
$\sigma\in\mathcal{B}(\Delta,\gamma)$ and $0\leq{j}<\mu_\sigma$, 
and compute their associated rational function (\ref{eq28}). 
The counting function $h(y;z)$ is then obtained as follows. 

\vspace{9pt}\noindent
{\bf Input:} $y\in\Z^m\cap\overline{\gamma},\,z\in\C^n$.\\
{\bf Output} $\rho=h(y;z)$.\\
Set $\rho:=0$. For every $\sigma\in\mathcal{B}(\Delta,\gamma)$ :

$\bullet$ Compute 
$\xi[\jmath[y,\sigma],\sigma]:=y-A_{\sigma}
\lfloor{}A_\sigma^{-1}y\rfloor\in\Z^m$.

$\bullet$ Read the value 
$R_1(\xi[\jmath[y,\sigma],\sigma],\sigma;z)/R_2(\sigma;z)$, 
and update $\rho$ by:
$$\rho:=\rho+\frac{R_1(\xi[\jmath[y,\sigma],\sigma],\sigma;z)}
{R_2(\sigma;z)}\,z_\sigma^{\lfloor{}A_\sigma^{-1}y\rfloor}.$$

For the whole space $\Z^m$ it suffices to consider {\it all} chambers 
$\gamma$ and all cones $C(\om(\xi[j,\sigma]),\sigma)\subset\R^n$, 
where $\sigma\in\mathcal{B}(\Delta,\gamma)$ and $0\leq{j}<\mu_\sigma$.

Finally, in view of (\ref{eq11})-(\ref{eq12}), the above algorithm 
can be symbolic, i.e., $z\in\C^m$ can be treated symbolically, and 
$\rho$ becomes a rational fraction of $z$.

\section{Generating function}
\label{generating}

An appropriate tool for computing the exact value of $h(y;z)$ 
in~(\ref{eq2}) is the formal generating function $H:\C^m\to\C$,
\begin{equation}\label{eqn1}
s\mapsto{H}(s)\,:=\,\sum_{y\in\Z^m}\,h(y;z)\,s^y\,=\,\prod_{k=1}^n\frac{1}{1-z_ks^{A_k}},
\end{equation}
where $s^y$ is defined in~(\ref{eq8}) and the sum is understood 
as a formal power series, so that we need not consider 
conditions for convergence. This generating function was already 
considered in Brion and Vergne \cite{brion3} with 
$\lambda=\ln\langle{s}\rangle$.

Following notation of \cite[p. 805]{brion3},
let $0\leq\hat{x}\in\R^n$ be a \textit{regular} vector, i.e., no entry $[A_\sigma^{-1}A\hat{x}]_j$ 
vanishes for any basis $\sigma\in\J_A$ or index $1\leq{j}\leq{m}$. 
Define~:
\begin{equation}\label{eqn3}
\varepsilon_{j,\sigma}\,:=\,\left\{\begin{array}{cl}
1&\hbox{if}\quad[A_\sigma^{-1}A\hat{x}]_j>0,\\
-1&\hbox{if}\quad[A_\sigma^{-1}A\hat{x}]_j<0.
\end{array}\right.\end{equation}

Next, for every basis $\sigma\in\J_A$, index $j\in\sigma$ 
and vector $u_{\not\sigma}\in\Z^{n-m}$, fix~:
\begin{equation}\label{eqn4}
\begin{array}{l}
\theta[j,\sigma,u_{\not\sigma}]\in\Z:\quad\hbox{the smallest integer greater}\\
\hbox{than or equal to}\quad-\varepsilon_{j,\sigma}[A_\sigma^{-1}
A_{\not\sigma}u_{\not\sigma}]_j.
\end{array}\end{equation}

Define also the vector $\eta[\sigma,u_{\not\sigma}]\in\Z^n$ by~:
\begin{equation}\label{eqn-eta}
\eta[\sigma,u_{\not\sigma}]_j\,=\,\left\{\begin{array}{cl}
u_j&\hbox{if}\;j\not\in\sigma;\\
\theta[j,\sigma,u_{\not\sigma}]&\hbox{if}\;
j\in\sigma,\;\varepsilon_{j,\sigma}=1;\\
1-\theta[j,\sigma,u_{\not\sigma}]&\hbox{if}\;
j\in\sigma,\;\varepsilon_{j,\sigma}=-1.\\
\end{array}\right.\end{equation}

The following expansion can be deduced from \cite{brion3}.

\begin{theorem}\label{expansion}
Let $0\leq\hat{x}\in\R^n$ be \textit{regular}
and consider the vectors 
$\eta[\sigma,u_{\not\sigma}]\in\Z^n$ defined in~(\ref{eqn-eta}) 
for $\sigma\in\J_A$ and $u_{\not\sigma}\in\Z^{n-m}$. The following
expansion holds: 
\begin{eqnarray}\label{eqn5}
&&\prod_{k=1}^n\frac{1}{1-z_ks^{A_k}}\;=\;\sum_{\sigma\in\J_A}
\Bigg[\prod_{j\in\sigma}\frac{1}{1-z_js^{A_j}}\Bigg]\times\\
\nonumber&&\qquad\times\,\frac{1}{R_2(\sigma;z)}\,
\sum_{u_{\not\sigma}\in\Z_{\mu_\sigma}^{n-m}}
z^{\eta[\sigma,u_{\not\sigma}]}\,s^{A\eta[\sigma,u_{\not\sigma}]},
\end{eqnarray}
where $\Z_{\mu_\sigma}=\{0,1,\ldots,\mu_\sigma-1\}$,
$\mu_\sigma=|\det{}A_\sigma|$ and:
\begin{equation}\label{eqn6}
z\,\mapsto\,R_2(\sigma;z)\,:=\,\prod_{k\notin\sigma}\left[1-
\bigl(z_kz_\sigma^{-A_\sigma^{-1}\!A_k}\bigr)^{\mu_\sigma}\right].
\end{equation}
\end{theorem}

\begin{proof}
From Brion and Vergne's identity
\cite[p.~813]{brion3}, 
\begin{equation}\label{key}
\prod_{j=1}^n\frac{1}{1-\e^{w_k}}\,=\,\sum_{\sigma\in\J_A}
\bigg[\prod_{j\in\sigma}\varepsilon_{j,\sigma}\bigg]\,
F(C^\sigma_{\hat{x}}+\rho_\sigma(C),L),
\end{equation}
where $F(C^\sigma_{\hat{x}}+\rho_\sigma(C),L)$ is the formal power 
series $\sum_l\e^l$ added over all elements $l$ in the intersection 
of the cone $C^\sigma_{\hat{x}}+\rho_\sigma(C)$ with the integer 
lattice $L=\Z[w_1,\ldots,w_n]$. Moreover, the coefficients 
$\varepsilon_{j,\sigma}$ are defined in~(\ref{eqn3}) and the 
cone $C_{\hat{x}}^\sigma$ is defined by the following formula 
\cite[p.~805]{brion3},
\begin{equation}\label{eqn7}
C^\sigma_{\hat{x}}\,=\,\bigg\{\sum_{j\in\sigma}\varepsilon_{j,\sigma}
\,x_j\,w_j\:\bigg|\:x_\sigma\in\R^m,\,x_\sigma\geq0\bigg\}.
\end{equation}

Finally, given the real vector space $W=\R[w_1,\ldots,w_n]$, 
every $\rho_\sigma:W\to{W}$ is a linear mapping defined by 
its action on each basis element $w_k$ of $W$,
\begin{equation}\label{eqn8}
\rho_\sigma(w_k)\,:=\,w_k\,-\,\sum_{j\in\sigma}\,[A_\sigma^{-1}\!A_k]_jw_j.
\end{equation}

Hence, $\rho_\sigma(w_j)=0$ for every $j\in\sigma$, and the 
cone $\rho_\sigma(C)$ is given by 
\begin{equation}\label{eqn9}
\rho_\sigma(C)\,=\,\bigg\{\sum_{k{\not\in}\sigma}x_kw_k\,-\,
\sum_{j\in\sigma}[A_\sigma^{-1}\!A_{\not\sigma}x_{\not\sigma}]_jw_j
\:\bigg|\:\begin{array}{l}
x_{\not\sigma}\in\R^{n-m},\\
x_{\not\sigma}\geq0
\end{array}\bigg\};
\end{equation}
see \cite[p.805]{brion3}.
Thus, every element in the intersection of the cone 
$C_{\hat{x}}^\sigma+\rho_\sigma(C)$ with the lattice 
$\Z[w_1,\ldots,w_n]$ must be of the form~:
\begin{eqnarray}\label{eqn10}
&&\sum_{k\not\in\sigma}x_kw_k\,+\,
\sum_{j\in\sigma}\varepsilon_{j,\sigma}\,\xi_j\,w_j,
\quad\hbox{with}\quad x_{\not\sigma}\in\N^{n-m},\\
\nonumber&&\qquad\xi_\sigma\in\Z^m\quad\hbox{and}\quad\xi_j\geq-
\varepsilon_{j,\sigma}[A_\sigma^{-1}\!A_{\not\sigma}x_{\not\sigma}]_j.
\end{eqnarray}

On the other hand, for every basis $\sigma$, 
define $\mu_\sigma=|\det{A}_\sigma|$ and~:
\begin{equation}\label{eqn11}
x_{\not\sigma}\,=\,u_{\not\sigma}+\mu_{\sigma}v_{\not\sigma},
\quad\hbox{with}\quad{}u_{\not\sigma}\in\Z_{\mu_\sigma}^{n-m}
\quad\hbox{and}\quad{}v_{\not\sigma}\in\N^{n-m}.
\end{equation} 

Moreover, as in (\ref{eqn4}), fix 
$\theta[j,\sigma,u_{\not\sigma}]\in\Z$ to be 
the smallest integer greater than or equal to 
$-\varepsilon_{j,\sigma}[A_\sigma^{-1}A_{\not\sigma}u_{\not\sigma}]_j$.
Thus, we can rewrite (\ref{eqn10}) so that the intersection 
of the cone $C_{\hat{x}}^\sigma+\rho_\sigma(C)$ with the lattice $\Z[w_1,\ldots,w_n]$ must be of the form~:
\begin{eqnarray}\label{eqn12}
&&\sum_{k\not\in\sigma}\big[u_kw_k+v_k\mu_\sigma\rho(w_k)\big]
\,+\,\sum_{j\in\sigma}\varepsilon_{j,\sigma}w_j
\big[\theta[j,\sigma,u_{\not\sigma}]_j+q_j\big],\\
\nonumber&&\qquad\hbox{with}\quad{}u_{\not\sigma}\in\Z_{\mu_\sigma}^{n-m},
\quad{}v_{\not\sigma}\in\N^{n-m}\quad\hbox{and}\quad{}q_\sigma\in\N^m.
\end{eqnarray}

We can deduce (\ref{eqn12}) from (\ref{eqn10}) by recalling the 
definition (\ref{eqn8}) of $\rho_\sigma(w_k)$ and letting~:
$$\xi_j\,:=\,\theta[j,\sigma,u_{\not\sigma}]\,+\,q_j\,-\,
\varepsilon_{j,\sigma}\mu_\sigma[A_\sigma^{-1}\!A_{\not\sigma}v_{\not\sigma}]_j.$$

Since $F(C^\sigma_{\hat{x}}+\rho_\sigma(C),L)$ is the formal power 
series $\sum_l\e^l$ with summation over all elements $l$ in~(\ref{eqn12}),
one obtains
\begin{eqnarray}\label{eqn13}
&&F(C^\sigma_{\hat{x}}+\rho_\sigma(C),L)\;=\\
\nonumber&&
\sum_{u_{\not\sigma}\in\Z_{\mu_\sigma}^{n-m}}\Bigg[\prod_{j\in\sigma}
\frac{\e^{\varepsilon_{j,\sigma}\theta[j,\sigma,u_{\not\sigma}]w_j}}
{1-\e^{\varepsilon_{j,\sigma}w_j}}\Bigg]\Bigg[\prod_{k\not\in\sigma}
\frac{\e^{u_kw_k}}{1-\e^{\mu_\sigma\rho_\sigma(w_k)}}\Bigg].
\end{eqnarray}

With $\eta[\sigma,u_{\not\sigma}]\in\Z^n$ being as in (\ref{eqn-eta}),
using (\ref{eqn13}) into (\ref{key}) yields the expansion
\begin{eqnarray}\label{eqn14}
&&\prod_{j=1}^n\frac{1}{1-\e^{w_k}}\;=\;\sum_{\sigma\in\J_A}
\;\sum_{u_{\not\sigma}\in\Z_{\mu_\sigma}^{n-m}}
\Bigg[\prod_{j\in\sigma}\frac{1}{1-\e^{w_j}}\Bigg]\times\\
\nonumber&&\qquad\times\Bigg[\prod_{j=1}^n
\e^{\eta[\sigma,u_{\not\sigma}]_jw_j}\Bigg]\Bigg[\prod_{k\not\in\sigma}
\frac{1}{1-\e^{\mu_\sigma\rho_\sigma(w_k)}}\Bigg].
\end{eqnarray}

Finally, we defined $w_k:=\ln(z_k)+\ln\langle{s}\rangle{}A_k$ for 
every index $1\leq{k}\leq{n}$, where the vectors $s,z\in\C^n$ have 
all their entries different from zero and $\ln\langle{s}\rangle$ is 
the $[1\times{n}]$ matrix defined in (\ref{eq8}). So
$\e^{w_k}=z_ks^{A_k}$. Moreover, recalling the definition (\ref{eqn8}) of 
$\rho_\sigma(w_k)$, the following
identities hold for all $1\leq{k}\leq{n}$,
\begin{equation}\label{eqn15}
\rho_\sigma(w_k)\,=\,\ln(z_k)\,-\,
\sum_{j\in\sigma}\ln(z_j)[A_\sigma^{-1}\!A_k]_j.
\end{equation}

Notice $\sum_{j\in\sigma}A_j[A_\sigma^{-1}A_k]_j=A_k$. A direct application of (\ref{eqn15}) and 
the identities $\e^{w_k}=z_ks^{A_k}$ yields (\ref{eqn5}), i.e.:
$$\prod_{k=1}^n\frac{1}{1-z_ks^{A_k}}=
\sum_{\sigma\in\J_A}\,\sum_{u_{\not\sigma}\in\Z_{\mu_\sigma}^{n-m}}
\frac{z^{\eta[\sigma,u_{\not\sigma}]}s^{A\eta[\sigma,u_{\not\sigma}]}}
{R_2(\sigma;z)}\prod_{j\in\sigma}\frac{1}{1-z_js^{A_j}},$$
$$with \quad R_2(\sigma;z)\,=\,\prod_{k\notin\sigma}\left[1-
\bigl(z_kz_\sigma^{-A_\sigma^{-1}\!A_k}\bigr)^{\mu_\sigma}\right].$$
\end{proof}
A direct expansion of (\ref{eqn5}) 
yields the following:

\begin{theorem}\label{general-counting}
Let $0\leq\hat{x}\in\R^n$ be \textit{regular}, and let $h$ and $\eta$ be as in (\ref{eq2}) and (\ref{eqn-eta}), 
respectively. Let $\J_A$ be the set of bases associated with $A$.
Then for every pair of $(y,z)\in\Z^m\times\C^n$ 
with $\|z\|<1$:
\begin{eqnarray}\label{way2}
h(y;z)&=&\sum_{\sigma\in\J_A,\,A_\sigma^{-1}y\geq0}
\frac{z_\sigma^{A_\sigma^{-1}y}}{R_2(\sigma;z)}
\sum_{u\in\Z_{\mu_\sigma}^{n-m}}\frac{z_{\not\sigma}^u}
{z_\sigma^{A_\sigma^{-1}\!A_{\not\sigma}u}}\times\\
\nonumber&&\quad\times\left\{\begin{array}{cl}1&\hbox{if}\;
A_\sigma^{-1}\big(y-A\eta[\sigma,u]\big)\in\N^m,\\
0&\hbox{otherwise,}\end{array}\right.
\end{eqnarray}
where: \hspace{1.8cm}$\Z_{\mu_\sigma}=\{0,1,\ldots,\mu_\sigma-1\}$,
$\mu_\sigma=|\det{}A_\sigma|$, 
\begin{eqnarray}\label{eqn20}
&&0\,\leq\,\big[A_\sigma^{-1}A\eta[\sigma,u]\big]_j
\,\leq\,1\quad\hbox{for each}\quad{}j\in\sigma\\
&&\hbox{and}\quad{}R_2(\sigma;z)\,:=\,\prod_{k\notin\sigma}\left[1-
\bigl(z_kz_\sigma^{-A_\sigma^{-1}\!A_k}\bigr)^{\mu_\sigma}\right].
\end{eqnarray}
\end{theorem}
\begin{proof}
Recall the expansion of $H(s)$ as a formal power series~:
\begin{equation}\label{eqn21}
H(s)\,=\,\sum_{y\in\Z^m}h(y;z)\,s^y\,=\,
\prod_{k=1}^n\frac{1}{1-z_ks^{A_k}}.
\end{equation}

We also have a similar formal power series for the product~:
$$\prod_{j\in\sigma}\frac{1}{1-z_js^{A_j}}\,=\prod_{j\in\sigma}
\Bigg[\sum_{q_j=1}^{\infty}z_j^{q_j}s^{A_jq_j}\Bigg]=
\sum_{q_\sigma\in\N^m}z_\sigma^{q_\sigma}s^{A_{\sigma}q_\sigma}.$$
Combining the latter with (\ref{eqn5}) yields:
\begin{equation}\label{eqn22}
H(s)\,=\,\sum_{\sigma\in\J_A}\;\sum_{q_\sigma\in\N^m}
\;\sum_{u_{\not\sigma}\in\Z_{\mu_\sigma}^{n-m}}
\frac{z^{\eta[\sigma,u_{\not\sigma}]}z_\sigma^{q_\sigma}
s^{A\eta[\sigma,u_{\not\sigma}]}s^{A_{\sigma}q_\sigma}}
{R_2(\sigma;z)}.
\end{equation}

Notice that (\ref{eqn21}) and (\ref{eqn22}) are identical.
Hence, if we want to obtain the exact value of $h(y;z)$ 
from (\ref{eqn22}), we only have to sum up all terms with exponent
$A\eta[\sigma,u_{\not\sigma}]+A_{\sigma}q_{\sigma}$ equal to $y$.
That is, recalling that each $A_\sigma$ is invertible,
\begin{eqnarray}\label{eqn24}
h(y;z)&=&\sum_{\sigma\in\J_A}\;
\sum_{u_{\not\sigma}\in\Z_{\mu_\sigma}^{n-m}}
\frac{z^{\eta[\sigma,u_{\not\sigma}]}z_\sigma^{q_\sigma}}
{R_2(\sigma;z)}\times\\
\nonumber&&\times\left\{\begin{array}{cl}1&\hbox{if}\;
A_\sigma^{-1}\big(y-A\eta[\sigma,u_{\not\sigma}]\big)\in\N^m,\\
0&\hbox{otherwise.}\end{array}\right.
\end{eqnarray}

On the other hand, setting 
$q_\sigma:=A_\sigma^{-1}(y-A\eta[\sigma,u_{\not\sigma}])$ 
and recalling the definition (\ref{eqn-eta}), 
\begin{equation}\label{eqn25}
z^{\eta[\sigma,u_{\not\sigma}]}z_\sigma^{q_\sigma}\,=\,
\frac{z_\sigma^{\eta[\sigma,u_{\not\sigma}]_\sigma}
z_{\not\sigma}^{u_{\not\sigma}}z_\sigma^{A_\sigma^{-1}y}}
{z_\sigma^{A_\sigma^{-1}\!A\eta[\sigma,u_{\not\sigma}]}}\,=\,
\frac{z_\sigma^{A_\sigma^{-1}y}z_{\not\sigma}^{u_{\not\sigma}}}
{z_\sigma^{A_\sigma^{-1}\!A_{\not\sigma}u_{\not\sigma}}}.
\end{equation}

We finally prove that the vector 
$A_\sigma^{-1}A\eta[\sigma,u_{\not\sigma}]$ 
is bounded, so that (\ref{eqn20}) holds,
\begin{equation}\label{eqn26}
0\,\leq\,\big[A_\sigma^{-1}A\eta[\sigma,u_{\not\sigma}]\big]_j
\,\leq\,1\quad\hbox{for each}\quad{}j\in\sigma.
\end{equation}

From (\ref{eqn4}) and (\ref{eqn-eta}),
the following identity
\begin{equation}\label{eqn27}
\big[A_\sigma^{-1}\!A\eta[\sigma,u_{\not\sigma}]\big]_j=
\big[A_\sigma^{-1}\!A_{\not\sigma}u_{\not\sigma}\big]_j+
\left\{\begin{array}{cl}
\theta[j,\sigma,u_{\not\sigma}]&\hbox{if}\;\varepsilon_{j,\sigma}=1;\\
1-\theta[j,\sigma,u_{\not\sigma}]&\hbox{if}\;\varepsilon_{j,\sigma}=-1
\end{array}\right.
\end{equation}
holds for every index $j\in\sigma$.

Next, suppose that $[A_\sigma^{-1}A_{\not\sigma}u_{\not\sigma}]_j=R+r$ where 
$R\in\Z$ and $0\leq{r}<1$ are the respective integer and fractional 
parts. We can obtain $\theta[j,\sigma,u_{\not\sigma}]$ in 
(\ref{eqn4}) as follows: If $\varepsilon_{j,\sigma}=1$ then
$\theta[j,\sigma,u_{\not\sigma}]=-R$ is the smallest integer 
greater than or equal to $-R-r$, and so (\ref{eqn27}) yields
$$\big[A_\sigma^{-1}\!A\eta[\sigma,u_{\not\sigma}]\big]_j\,=\,r.$$

If $\varepsilon_{j,\sigma}=-1$ and $0<r<1$, then
$\theta[j,\sigma,u_{\not\sigma}]=R+1$ and~:
$$\big[A_\sigma^{-1}\!A\eta[\sigma,u_{\not\sigma}]\big]_j\,=\,r.$$

At last, if $\varepsilon_{j,\sigma}=-1$ and $r=0$, then
$\theta[j,\sigma,u_{\not\sigma}]=R$ and so,
$$\big[A_\sigma^{-1}\!A\eta[\sigma,u_{\not\sigma}]\big]_j\,=\,1.$$

Therefore (\ref{eqn26}) holds because $0\leq{r}<1$. Finally,
a direct application of (\ref{eqn25}) into (\ref{eqn24}) yields the 
following version of (\ref{way2}): 
\begin{eqnarray*}
h(y;z)&=&\sum_{\sigma\in\J_A,\,A_\sigma^{-1}y\geq0}\frac{1}
{R_2(\sigma;z)}\sum_{u_{\not\sigma}\in\Z_{\mu_\sigma}^{n-m}}
\frac{z_\sigma^{A_\sigma^{-1}y}z_{\not\sigma}^{u_{\not\sigma}}}
{z_\sigma^{A_\sigma^{-1}\!A_{\not\sigma}u_{\not\sigma}}}\times\\
&&\qquad\times\left\{\begin{array}{cl}1&\hbox{if}\;
A_\sigma^{-1}\big(y-A\eta[\sigma,u_{\not\sigma}]\big)\in\N^m,\\
0&\hbox{otherwise.}\end{array}\right.
\end{eqnarray*}
Notice that the first sum is 
calculated only on those bases $\sigma\in\J_A$ which satisfy 
$A_\sigma^{-1}y\geq0$. This follows from (\ref{eqn24}) combined
with (\ref{eqn26}).
\end{proof}

Observe that (\ref{way2}) is different from (\ref{eq15}) or
(\ref{eq41}) because in (\ref{eq15}) and (\ref{eq41}), the
summation is over bases $\sigma$ in the subset
$\mathcal{B}(\Delta,\gamma)\subset\{\J_A;\,A_\sigma^{-1}y\geq0\}$.

%
%

\end{document}